\documentclass[11pt]{amsart} 

\usepackage{euscript}
\usepackage{amssymb,amscd}

      \theoremstyle{plain}
      \newtheorem{theorem}{Theorem}[section]
      \newtheorem{lemma}[theorem]{Lemma}
      \newtheorem{corollary}[theorem]{Corollary}
      \newtheorem{proposition}[theorem]{Proposition}

      \makeatletter
      \def\@setcopyright{}
      \def\serieslogo@{}
      \makeatother

\newcommand{\R}{\mathbb R}
\newcommand{\Rk}{\mathbb R^k}
\newcommand{\C}{\mathbb C}
\newcommand{\Z}{\mathbb Z}
\newcommand{\Zk}{\mathbb Z^k}
\newcommand{\Q}{\mathbb Q}
\newcommand{\T}{\mathbb T}
\newcommand{\m}{\mathbb \mu ^F}
\newcommand{\Tm}{\mathbb T^m}

\def \a{\alpha}
\def \G{\Gamma}
\def \g{\mathfrak{g}}
\def \h{\mathfrak{h}}
\def \L{\Lambda}
\def \l{\lambda}

\def \A{(\EuScript{A})}
\def \M{\EuScript{M}}

\begin{document}

\author[Boris Kalinin   and  Ralf Spatzier]
{Boris Kalinin  and Ralf Spatzier }

  \address{The University of Michigan, Ann Arbor, MI}

  \email{kalinin@umich.edu, spatzier@umich.edu }

\title[Measurable rigidity for higher rank abelian actions]
 {Rigidity of the measurable structure for algebraic actions of higher rank abelian groups}
\thanks{The first author was partially supported by  NSF grant DMS 0140513, the second by 
NSF grants DMS 0203735 and 9971556}

\begin{abstract}

We investigate rigidity of measurable structure for higher rank abelian algebraic actions. In particular, 
we show that ergodic measures for these actions fiber over a 0 entropy measure
with Haar measures along the leaves. We deduce various rigidity theorems for 
isomorphisms and joinings as corollaries. 

\end{abstract}

\maketitle

\section{Introduction}

We consider algebraic  actions of higher rank abelian groups on 
homogeneous spaces. We establish rigidity properties of the measurable structure for
these actions. This is part of a general program investigating such actions and their
rigid structures such as local rigidity and  triviality of cocycles \cite{KS1,KS2}. As for
the measurable structure for Anosov actions without rank one factors, 
one conjectures that invariant ergodic measures for such
actions are always algebraic, i.e. Haar measures on closed homogeneous subspaces
\cite{Margulis,KS3}. For the case of measures on the circle 
invariant under $x2,x3$ this is the famous Furstenburg conjecture \cite{F}.
To this date, the only results known in this
direction concern  the classification of invariant measures with positive entropy. 
The first such results      were achieved 
by  Lyons and later Rudolph and Johnsson for the Furstenburg conjecture
\cite{Lyons,Rudolph,J1,J2}. More general results were obtained
later by Katok and Spatzier in \cite{KS3,KS4} using geometric ideas. 

In this paper we generalize the approach of \cite{KS3} to  general partially
hyperbolic algebraic actions assuming that certain one-parameter subgroups have
``large'' ergodic components. As a  main application we show that measurable
isomorphisms between two algebraic  actions of this type have to be algebraic. This had so far only been
proved for the case of a partially hyperbolic $\Zk$-action  by toral automorphisms
\cite{KKS}. 
Our result covers  actions on nilmanifolds by automorphisms as well as
Weyl chamber flows and other 
actions on homogeneous spaces of higher rank semisimple groups. We also obtain results
on the rigidity of  joinings. In particular, joinings in the semisimple case have to be
algebraic assuming suitable one-parameter subgroups are ergodic. 
Previously,   rigidity results for joinings  had only been obtained for actions
by toral automorphisms \cite{KaK}. Even in that case, our results apply more generally and
resolve some questions raised in  \cite{KaK}.

We state the main results for $\Rk$-actions in Section 3  and for $\Zk$-actions in Section 4 after 
establishing  basic notions in Section 2. The proof in the case of an $\Rk$-action on a homogeneous
space $G/\Gamma$ is given in Sections 5 and 6. The central idea  of the proof is the same as in \cite{KS3}. 
Essentially we prove a dichotomy for the conditional measures of an invariant measure $\mu$ along various coarse Lyapunov foliations.
They are either Haar or atomic. This is done as in \cite{KS3} by moving these conditional measures by
suitable elements in $\Rk$ which act isometrically along these foliations. Now either all these conditional
measures are atomic and the entropy is 0 or the measure is invariant under unipotent groups. 
 Ratner's theorem then shows that suitable ergodic components of the measure for the unipotent action 
 are Haar measures on  compact orbits of some subgroup $H$ of $G$. We now show that $\mu$ 
 essentially lives on an orbit or orbit closure of the normalizer of
 $H$. From this we can describe $\mu$ as a measure on a fibration with $H$-orbits and Haar measures as
 fibers.  In the semisimple and nilpotent case, we can additionally show that this quotient is algebraic. 
  The main novelty of our argument lies in this latter part of the proof. The factorization obtained is 
 quite strong and allows us in particular to get general results about isomorphism rigidity and joinings. 
Our techniques give new results about joinings even for toral automorphisms.
  We  derive the  results for measurable isomorphisms and joinings  as well as the 
 results for $\Zk$-actions in the later sections of the paper.

  We thank G. Prasad and D. Witte for various helpful discussions,  and the latter especially for the proof of
  Proposition 8.1.

\section{Preliminaries and Notations.}

\subsection{Algebraic $\Rk$ Actions.}
Let $G$ be a simply connected Lie group and $\Gamma \subset G$ be a uniform lattice.
Let $\rho$ be an embedding of $\Rk$ into $G$ and let $A=\rho(\Rk) \subset G$ be its
image. Then $\Rk$ acts on $\M=G/\G$ via left translations by corresponding elements of
$A$. We denote this action by $\a$. The linear part of $\a$ is the adjoint 
representation of $A$ on the Lie algebra ${\g}$ of $G$. We will assume that for all 
$a \in A$, $Ad(a)$ is a diagonalizable (semisimple) automorphism of the Lie algebra 
$\g$. We note that $\alpha$ preserves the Haar measure  $\lambda$ on $\M=G/\G$.
We call $(\alpha, \mu)$ an algebraic $\Rk$ action. A prime example of such an action
is given by an embedding of $\Rk$ into a split Cartan subgroup of a semisimple Lie
group $G$, for example $SL(n,\R)$ with $n > k$.

Let us briefly describe the relation between the Lyapunov exponents and the algebraic
structure of the action.
Since $A$ is commutative the Lie algebra $\g$ splits into invariant subspaces 
for the adjoint action of $A$, which are the eigenspaces corresponding to real 
eigenvalues and invariant subspaces corresponding to the pairs of complex eigenvalues.
For each invariant subspace let us denote by $\chi(a)$ the logarithm of the modulus 
of the corresponding eigenvalue of $Ad(a)$. Then $\chi(a)$ is a liner functional
on $A$. This splitting gives rise to the right-invariant splitting of the tangent bundle 
of $G$, which projects to the $\a$ invariant splitting of the tangent bundle of 
$\M=G/\G$. The latter is a refinement of the Lyapunov decomposition of $T\M$ into 
Lyapunov subspaces for the action $\a$. Indeed, it is easy to see that for any
vector $v$ in one of the right-invariant distributions its Lyapunov exponent 
$\chi (\a (a),v)$ is equal to the corresponding $\chi(a)$.

In fact it is easy to see that we can take an inner product on $\g$ such that
the invariant subspaces are orthogonal and for any $a \in A$ the restriction of 
$Ad(a)$ to each of the invariant subspaces is a scalar multiple of isometry.
This inner product gives rise to the right-invariant Riemannian metric on  
$G$, which projects to the Riemannian metric on $\M=G/\G$. It follows
that with respect to this metric 
$$||D\a(a)v|| = e^{\chi(a)}||v||$$
for any $a \in \Rk$ and vector $v$ in the distribution corresponding to $\chi$.

Each Lyapunov exponent for the action $\a$ is a linear functional on $A$. 
Its positive half--space is called  a {\it Lyapunov half--space} and its boundary 
a {\it Lyapunov hyperplane}. Note that Lyapunov exponents  may be proportional to 
each other with positive or negative coefficients. In this  case, they define the 
same Lyapunov hyperplane. Combining Lyapunov distributions  with the same
Lyapunov half--space we obtain {\it coarse Lyapunov distributions}
which form the {\it coarse Lyapunov decomposition}.

An element $a \in \Rk$ is called {\it regular} 
if it does not belong to any Lyapunov hyperplane. All other elements are 
called {\it singular}. Call a singular element {\it generic} if it  belongs 
to only one Lyapunov hyperplane.

For an element $a \in A$  the {\em stable, unstable} and  {\em
neutral\/} distributions $E^- _a, E^+ _a$ and $E^0 _a$  are defined as
the sum of the Lyapunov  spaces for which the value of the
corresponding Lyapunov exponent on $a$ is  negative, positive and 0
respectively. From the above discussion we see that
the derivative of any singular element $a$ acts
isometrically on its neutral distribution.  So in this case $E^0 _a$
coincides with the isometric distribution $E^I _a$.

Notice that while some Lyapunov distributions may not be integrable
the stable and unstable distributions always integrate to homogeneous 
foliation whose leaves are in fact cosets of nilpotent subgroups.
For an element $a$ we will denote the integral foliations of the stable, 
and unstable distributions $E^- _a$ and $E^+ _a$ by $W^-_a$ and $W^+ _a$.

We note that $E^- _a$ is the sum of certain coarse Lyapunov distributions.  
Conversely, any coarse Lyapunov distribution is the intersection of the
strong stable distributions for certain elements of the action. Thus, any 
coarse Lyapunov distribution integrates to a homogeneous foliation whose
leaves are also cosets of a nilpotent subgroup.

\subsection{Algebraic $\Zk$ Actions and Their Suspensions.}
Let $G$ be a simply connected nilpotent Lie group and $\Gamma \subset G$ be 
a uniform lattice. If $A$ is an automorphism of $A$ such that $A(\G)=\G$, than 
it projects to an automorphism $F_A$ of the nilmanifold $\M=G/\G$ which
preserves the Haar measure $\l$ on $\M$. The differential of $A$ at $e \in G$
is an automorphism of the Lie algebra $\g$ of $G$ and is called the linear
part of $A$. By a theorem of W. Parry (\cite{P}) $F_A$ is ergodic with respect 
to $\l$ if and only if no eigenvalue of the quotient of the linear part of $A$
on the abelianization $\g/[\g,\g]$ is a root of unity, furthermore, in this 
case $F_A$ is a K automorphism. 

If automorphisms $A_1$, ..., $A_k$ as above commute, they define a $\Zk$ action
$\a$ by automorphisms $F_{A_i}$ of $G/\G$. Such an action we call an algebraic
$\Zk$ action by automorphisms of the nilmanifold $\M=G/\G$.  
A prime example of such an action is given by an embedding of $\Zk$ into 
$SL(n,\Z)$, with $n > k$, whose elements act by automorphisms of $\R^n$
which preserve the lattice $\Z^n$ and thus project to the automorphisms
of the torus $\T^n$. In this paper we will restrict our attention to the 
actions for which the linear parts of the automorphisms are diagonalizable.

Let us briefly describe the Lyapunov exponents for the algebraic
$\Zk$ action by automorphisms of the nilmanifolds. Since $A_1$, ..., $A_k$ 
commute the Lie algebra $\g$ splits into invariant subspaces for their linear 
parts, which are the eigenspaces corresponding to real eigenvalues and invariant 
subspaces corresponding to the pairs of complex eigenvalues. 
For each invariant subspace let us denote by $\chi(a)$ the logarithm of the modulus 
of the corresponding eigenvalue of $Ad(a)$. Then $\chi(a)$ is an additive functional
on $\Zk$. This splitting gives rise to the right-invariant splitting of the tangent 
bundle of $G$, which projects to the $\a$ invariant splitting of the tangent bundle 
of $\M=G/\G$. The latter is again a refinement of the Lyapunov decomposition of 
$T\M$ into Lyapunov subspaces for the action $\a$. Similarly to the case of algebraic 
$\Rk$ actions, one can easily construct a Riemannian metric on $\M=G/\G$ with respect 
to which 
$$||D\a(a)v|| = e^{\chi(a)}||v||$$
for any $a\in \Zk$ and vector $v$ in the distribution corresponding to $\chi$.

We would like to view the additive functionals on $\Zk$ as linear functionals on $\Rk$
and be able to operate with the elements in their kernels, i.e. in the Lyapunov
hyperplanes. For this we need to pass from an action of $\Z^k$ to the corresponding 
action of $\R^k$ via the so-called suspension construction.

Suppose $\Z ^k$ acts on $G/\G$.  Embed $\Z ^k$ as a lattice in $\R^k$.
Let $\Z ^k$ act on $\R^k \times (G/\G)$  by $z\, (x,m) =(x-z,z \cdot m)$
and form the quotient \[ N = (\R^k \times (G/\G)) / \Z ^k.\] Note that the
action of $\R^k$  on $\R^k \times (G/\G)$ by $x \cdot(y,n)=(x+y,n)$
commutes with the $\Z ^k$-action and therefore descends to the
solvemanifold $N$. 

The manifold $N$ is a fibration over the "time" torus $\T^k$ with the
fiber $G/\G$.  We note that $TN$ splits into the direct sum $TN=T_fN
\oplus T_oN$ where $T_fN$ is the subbundle tangent to the $G/\G$ fibers
and $T_oN$ is the subbundle tangent to the orbit foliation.  The
Lyapunov exponent corresponding to $T_oN$ is always identically zero.
To exclude this trivial case, when we will consider only the Lyapunov 
exponents corresponding to $T_fN$. These Lyapunov exponents of the $\R^k$ 
action are the extensions of the Lyapunov exponents of the $\Z^k$ action 
to the linear functionals on $\R^k$. 

We note that the suspension of an algebraic $\Z ^k$-action is an algebraic 
$\Rk$ action. The further definitions and basic properties related to the
Lyapunov exponents, distributions and foliations are the same as in the
case of algebraic $\Rk$ actions. We call an action {\em totally nonsymplectic}
or TNS if no Lyapunov exponents are proportional with negative coefficients.

Note that any $\Z ^k$-invariant measure on $G/\G$ lifts to a unique
$\R^k$-invariant measure on $N$ and  conversely any invariant measure
for the suspension induces a unique  invariant measure  for the
original action.

\subsection{Measure--Theoretic Notations.}
Let $\a$ be a $\Zk$ or $\Rk$ action on a compact metric space. For any 
$\alpha$-invariant Borel probability  measure  $\mu$ we denote by  $\xi_a$  
the partition into ergodic components of the element $\a(a)$ and by $h_{\mu}(a)$ 
the measure--theoretic entropy of the map $\a(a)$. If $W$ is a continuous
foliation of the metric space, $\xi (W)$ denotes its measurable hull, i.e. the
measurable partition whose elements consist $\mu$-a.e. of complete leaves
of $W$.

Suppose $\a$ and $\a'$ are measurable actions of the same group $T$ by 
measure--preserving transformations of the spaces $(X,\l)$ and $(Y,\nu)$,
respectively. A probability measure $\mu$ on $X \times Y$ is called a 
{\it joining} of $\a$ and $\a'$ if it is invariant under the diagonal 
(product) action $\a \times \a'$ of $T$ on $X \times Y$ and its projections to 
$X$ and $Y$ coincide with $\l$ and $\nu$, respectively. The product measure
$\mu =\l \times \nu$ is called the trivial joining, and all other joining
measures are called nontrivial joinings.

A particular case of a joining is given by a measurable isomorphism, i.e. a measurable map
$\phi$ which intertwines the actions $\a$ and $\a'$. In fact, the push-forward of $\lambda$
under the graph of $\phi$, $(id \times \phi) _* (\lambda)$, is a joining. We call a 
measurable isomorphism
$\phi$ between two algebraic $\Zk$-actions $\a$ and $\a'$ {\em algebraic} if $\phi$ 
 lifts to an affine isomorphism $\psi$ between the ambient groups, i.e. $\psi$ is a composite
 of an isomorphism with left translations. 

For a probability measure $\mu$ on $G/\Gamma$ we denote by $\Lambda (\mu)$ 
the subgroup of all elements in $G$ whose left translations preserve $\mu$ and 
by $L(\mu)$ the Lie algebra of $\Lambda (\mu)$. The measure $\mu$ is called 
{\it algebraic} if there exists $x \in G/\Gamma$ such that $\mu (x \Lambda (\mu))=1$, 
i.e. $\mu$ is the Haar measure on a single left coset of $\Lambda (\mu)$.

%%%%%%%%%%%%%%%%%%%%%%%%%%%%%% MAIN RESULTS %%%%%%%%%%%%%%%%%%%%%%%%%%%%%%%%%%%

\section{Main Results.} \label{main results}

In this section we formulate our main results on rigidity of invariant measures,
measure--preserving isomorphisms, and joinings for homogeneous $\Rk$ actions.

\subsection{Rigidity of Invariant Measures}

\begin{theorem}  \label{im} 
Let $G$ be a  connected Lie group, $\Gamma \subset G$ be a uniform lattice,  
and $A$ be a subgroup of $G$ isomorphic to $\Rk$, with $k\geq 2$, such that for all 
$a \in A$, $Ad(a)$ is a diagonalizable over $\C$ automorphism of the Lie algebra $\g$. 
Let $\alpha$ denote the $A$-action on $G/\G$ by left translations, and let $\mu$ 
be an ergodic invariant measure for $\alpha$. 
Suppose further that either the eigenvalues of  $Ad(a)$ are real for all $a \in A$
or $\mu$ is weakly mixing for $\alpha$.

Suppose for some regular element $b \in A$ we have that for each coarse Lyapunov
subfoliation $W^P \subset W^- _b$ there exists a singular element $a$ in the
corresponding Lyapunov hyperplane $\partial P$ which act ergodically on $W^P$, i.e.
$\xi _{a} \leq \xi (W^P)$.

Then there exists a subgroup $H \subset G$ such that: 
 
\begin{enumerate}
\item $A$ normalizes $H$, i.e. $aHa^{-1}=H$ for any $a \in A$.
\item $H \subset \Lambda (\mu)$, i.e. left translations by elements of $H$
      leave $\mu$ invariant.
\item For almost every leaf of $W_b^-$ the conditional measure of $\mu$ 
      on this leaf is Haar on a single left coset of the Lie subgroup of $G$ 
      with Lie algebra $E_b^- \cap \h$, where $\h$ is the 
      Lie algebra of $H$. 
\item For $\mu$ a.e. $x$ the coset $Hx$ is compact in M, and the ergodic components
      of $\mu$ with respect to the left action of $H$ are Haar measures on the left
      cosets of $H$. 
      Moreover,  $\a$ induces a measurable factor action of $A$ on the space of
      these ergodic components. 
      The entropy of the factor measure is zero with respect to any element in the 
      Weyl chamber of $b$.
      
\item Suppose, in addition, that the coset $N_G(H)y$ of the normalizer of $H$ in $G$ is 
compact for some $y$ whose $A$ orbit is dense in $supp(\mu)$. Then
$supp(\mu) \subset N_G(H)y$ and for any $x \in N_G(H)y$ the coset $Hx$ 
      is compact in $G/\G$.
      In this case,  
      the restriction of $\a$ to $N_G(H)y$ induces an algebraic factor action  
      which is algebraically isomorphic to the action of 
      $y^{-1}Ay$ on $(y^{-1}Hy) \setminus y^{-1}N_G(H)y/ \G$. 
      The entropy of the factor measure is zero with respect to {\it any} element 
      in $A$ and the conditional measures in the fibers are Haar. Also, conclusion
      3. holds for the stable foliation $W_c^-$ of {\bf any} element $c \in A$.

\end{enumerate}
\end{theorem}
\smallskip
We prove this theorem in Sections \ref{proof im} and \ref{proof A}.
\smallskip

\noindent {\bf Remark 1}. If $H$ is a normal subgroup then $N_G(H)y=Gy$ and 
the last conclusion of the theorem above applies. 
\smallskip

\noindent {\bf Remark 2}. It easy to see that the measure $\mu$ in the theorem is 
algebraic, i.e. is supported on a single coset of $H$, if and only if $A \subset H$,
in other words if $H$ is $\a$ invariant. 
\smallskip

Theorem \ref{im} allows us to obtain some important corollaries,
which include the isomorphism rigidity theorem in Section \ref{Irig}.
\smallskip

\begin{corollary} \label{corK}
Suppose, in addition to the assumptions of Theorem \ref{im},
that for some element $c$ in the Weyl chamber of $b$, $\a (c)$ is a $K$-automorphism 
with respect to the measure $\mu$. Then $\mu$ is an algebraic measure, i.e. it is Haar on 
a single left coset $Hx$.  
\end{corollary}

\begin{proof}
Since $\a (c)$ is a K-automorphism with respect to the measure $\mu$, the measurable
factor with zero entropy has to be trivial. Hence $\mu$ is Haar on a single left 
coset $Hx$.
\end{proof}
\smallskip

\begin{corollary}[Semisimple Case] \label{corSSimple}
Let $G$ be a semisimple group and $\G \subset G$ be a uniform lattice. Let $\a$ 
be an $\Rk$-action, $k\geq 2$, on $G/\G$ by left translations given by an embedding 
$\rho$ of $\Rk$ into a split Cartan subgroup of $G$. Let $\mu$ be an ergodic invariant 
measure for $\a$.  

Suppose for some regular element $b \in \Rk$ we have that for each coarse Lyapunov
subfoliation $W^P \subset W^- _b$ there exists a singular element $a$ in the
corresponding Lyapunov hyperplane $\partial P$ which act ergodically on $W^P$, i.e.
$\xi _{a} \leq \xi (W^P).$ Then the conclusions (1)--(5) of Theorem \ref{im} hold.
\end{corollary}
\smallskip

We prove this corollary in Section \ref{proof corSSimple}. 
\smallskip

While the ergodicity assumption $\xi _{a} \leq \xi (W^P)$ in Theorem \ref{im}
and Corollary \ref{corSSimple}, and K-automorphism assumption in  Corollary \ref{corK}
are somewhat restrictive, they can be verified for a certain type of joining measure
associated with measure--preserving isomorphisms. This allows us to deduce the following
results on the rigidity of measure--preserving isomorphisms.
\smallskip

\subsection{Isomorphism Rigidity and Rigidity of Joinings}

We note that the ergodicity and K-automorphism assumptions in the following theorem 
are with respect to the Haar measure. Therefore they can be checked in various cases
as Corollary \ref{IrigSSimple} illustrates. 

\begin{theorem} \label{Irig}
Let $G$ be a  connected Lie group and $\Gamma \subset G$ be a uniform lattice.
Let $\rho$ be an embedding of $\Rk$ into $G$ such that for all $a \in \Rk$, 
$Ad(\rho(a))$ is a diagonalizable over $\C$ automorphism of the Lie algebra $\g$.
Let $\alpha$ denote the corresponding $\Rk$-action on $G/\G$ by left translations, 
and let  $\lambda$ be the Haar measure on $G/\G$.
Let $(\alpha', \lambda ')$ be an $\Rk$-action of the same type on $G'/\G'$.

Suppose that $\Rk$ contains a subgroup $S$ isomorphic to $\R ^2$ for which
\begin{enumerate}
\item The $\a$-action of any one-parameter subgroup is ergodic with respect 
      to $\lambda$.
\item There exists an element $b \in S$ such that $\a(b)$ is a K-automorphism 
      with respect to $\lambda$.
\end{enumerate}

Then any measure--preserving isomorphism $h: G/\Gamma \to G' /\Gamma '$  
between $(\a,\lambda)$ and $(\a', \lambda ')$ coincides \textup{(mod\;0)} 
with an affine isomorphism, and hence $\a$ and $\a '$ are algebraically 
isomorphic. 
\end{theorem}
\smallskip

\begin{proof} Consider the diagonal  action $\a \times \a'$ of $\R ^k$ 
on $G/\Gamma \times G' /\Gamma '$. Let $\mu$ be the push forward of $\lambda$ to
the graph of the measure--preserving isomorphism $h$. Then $\mu$ is an invariant 
measure for $\a \times \a'$, and $(\a \times \a ', \mu)$ has the same ergodic 
properties as $(\a,\lambda)$. Hence the assumptions of Theorem \ref{im} and 
Corollary \ref{corK} are satisfied for the restriction of $\a \times \a'$ to the
subgroup $S \subset \Rk$. Thus $\mu$ is Haar on a single left coset $Hx$. 
Since $h$ is a measure--preserving isomorphism and $\mu$ projects to the Haar 
measures in the factors, we see that the coset $Hx$ has to project one to one 
onto the factors. Thus $Hx$ is a graph of an affine isomorphism $\tilde h$ 
between $G/\Gamma$ and $G' /\Gamma '$ which coincides \textup{(mod\;0)} with $h$.
\end{proof}
\smallskip

\begin{corollary}[Semisimple Case] \label{IrigSSimple}
Let $G$ be a semisimple group without compact factors, $\G \subset G$ be an irreducible 
uniform lattice, and $\lambda$ be the Haar measure on $G/\G$. Let $(\a,\lambda)$ 
be an $\Rk$-action, $k\geq 2$, on $G/\G$ by left translations given by an 
embedding $\rho$ of $\Rk$ into a split Cartan subgroup of $G$. 

Let $G'$ be any simply connected Lie group and $\Gamma' \subset G'$ be a uniform 
lattice. Let $\rho'$ be an embedding of $\Rk$ into $G'$ such that for all $a \in \Rk$, 
$Ad(\rho'(a))$ is a diagonalizable over $\C$ automorphism of the Lie algebra $\g'$. 
Let $\alpha'$ denote the corresponding $\Rk$-action on $G'/\G'$ by left translations, 
and let $\lambda'$ be the Haar measure on $G'/\G'$.
Then any measure--preserving isomorphism $h: G/\Gamma \to G' / \Gamma '$ between 
$(\a,\lambda)$ and $(\a', \lambda ')$ coincides \textup{(mod\;0)} with an affine 
isomorphism, and hence $\a$ and $\a '$ are algebraically isomorphic.
\end{corollary}

\begin{proof}
We apply Theorem \ref{Irig}. It suffices to check that the restriction of $\a$ to any 
one-parameter subgroup of $\Rk$ is ergodic with respect to $\lambda$ and that 
there exists an element $b \in \Rk$ such that $\a (b)$ a K-automorphism with 
respect to $\lambda$. The first follows from the Howe-Moore theorem \cite{KS2} and the latter is proved  in 
\cite{D}.
\end{proof}
\medskip
 
In the proof of Theorem \ref{Irig} we applied the main Theorem \ref{im} to the
special type of joining measure. The following theorem describes the rigidity
of joinings in general. Here, however, we make the ergodicity assumption for the
joining measure itself.

\begin{theorem}[Semisimple Case] \label{JrigSSimple}
Let $G_1$ be a semisimple group, $\G_1 \subset G_1$ be 
a cocompact lattice, and $\lambda_1$ be the Haar measure on $G_1/\G_1$. 
Let($\a_1,\lambda_1)$ be an $\Rk$-action, $k\geq 2$, on $G_1/\G_1$ by left translations 
given by an embedding $\rho_1$ of $\Rk$ into a split Cartan subgroup of $G_1$.
Suppose that $\a_1(a)$ is a K-automorphism with respect to $\l_1$ for some $a \in \Rk$.

Let $(\alpha_2,\lambda_2)$ be an $\Rk$-action of the same type on $G_2/\G_2$, 
and let $(G_1/\Gamma_1 \times G_2 /\Gamma _2, \a_1 \times \a_2, \mu)$ 
be a nontrivial joining of the actions $\a_1$ and $\a_2$. 
If any one-parameter subgroup of $\Rk$ acts ergodically with respect to the joining 
measure $\mu$, then $\mu$ is algebraic and the actions $\alpha_1$ and $\alpha_2$ have 
algebraic factors whose finite covers are algebraically isomorphic. 
\end{theorem}

We prove this theorem in Section \ref{proof Jrig}.

%%%%%%%%%%%%%%%%%%%%%% Nilmanifold CASE %%%%%%%%%%%%%%%%%%%%%%%%%%%%%%%%%%%%%%%

\section{$\Zk$ Actions by Automorphisms of Nilmanifolds.}  \label{main results Nil}

In this section we state the counterparts of the rigidity results of Section 
\ref{main results} in the context of $\Zk$ actions by automorphisms of tori 
and nilmanifolds.
\smallskip

\subsection{Rigidity of Invariant Measures}  \label{measures Nil}

We will consider actions of $\Zk$ on nilmanifolds $G/ \Gamma$ by  automorphisms. If the automorphisms are diagonalizable, 
we can extend the homomorphism of $\Zk$ into the automorphisms of $G$ to $\Rk$. This allows us to consider the
induced action in the next theorem. 

\begin{theorem}  \label{imNil}
Let $\alpha$ be a $\; \Zk$, $k \ge 2$ action by automorphisms of $G/\Gamma$,
where $G$ is a simply connected nilpotent Lie group and $\G \subset G$ is a cocompact 
lattice. Let $\mu$ be an ergodic invariant measure for $\a$. Suppose that the
corresponding automorphisms of the Lie algebra $\g$ are diagonalizable over $\C$, 
and either all their eigenvalues are real or $\mu$ is a weakly mixing measure for 
$\alpha$.

Suppose for some regular element $b$ for the suspension action of $\Rk$
we have that for each coarse Lyapunov subfoliation $W^P \subset W^- _b$ 
there exists a singular element $a$ in the corresponding Lyapunov hyperplane 
$\partial P$ which acts ergodically on $W^P$, i.e.
$\xi _{a} \leq \xi (W^P).$

Then there exists a subgroup $H \subset G$ such that: 
 
\begin{enumerate}
\item $H$ is $\a$-invariant.
\item $H \subset \Lambda (\mu)$, i.e. left translations by elements of $H$
      leave $\mu$ invariant.
\item For any element $c \in \Zk$ and for almost every leaf of $W_c^-$ the 
      conditional measure of $\mu$ 
      on this leaf is Haar on a single left coset of the Lie subgroup of $G$ 
      with Lie algebra $E_c^- \cap \h$, where $\h$ is the 
      Lie algebra of $H$. 
\item $\mu$ is supported on a coset $N_G(H)y$ of the normalizer of $H$ in $G$. 
      For any $x \in N_G(H)y$ the coset $Hx$ is compact in $G/\G$, and the
      restriction of $\a$ to $N_G(H)y$ induces an algebraic factor action  
      which is isomorphic to $(y^{-1}Hy) \setminus y^{-1}N_G(H)y/ \G$. 
      The entropy of the factor measure is zero with respect to {\it any} 
      element in $\Zk$ 
      and the conditional measures in the fibers are Haar.
\end{enumerate}

\end{theorem}
\smallskip

This theorem is proved in Section \ref{proof imNil}. It generalizes the earlier
results for the actions by toral automorphisms, see Theorem 5.1' in \cite{KS4}
and its modified form Theorem 3.1 in \cite{KaK2}. In the toral case, however,
these theorems make no assumption on the corresponding automorphisms of the Lie 
algebra. While technically complicated, one can probably allow the presence of
Jordan blocks and thus remove the diagonalizability assumption in the nilmanifold
case. The presence of complex eigenvalues in the absence of weak mixing is
overcome in the toral case by, possibly, splitting the measure into finitely
many components. However, it is not clear whether the corresponding arguments 
in the proof of Theorem 3.1 in \cite{KaK2} can be generalized to the 
nonabelian case.

\medskip

In the case of actions by automorphisms of nilmanifolds the ergodicity assumption 
$\xi _{a} \leq \xi (W^P)$ in Theorem \ref{imNil} can be verified for an arbitrary
measure in the important special case of TNS actions, i.e. the partially hyperbolic
actions with no negatively proportional Lyapunov exponents. This is illustrated
by the following corollary.

\begin{corollary}[TNS case] \label{corTNS}
Let $\alpha$ be a TNS $\; \Zk$ action by automorphisms of $G/\Gamma$, where $G$ is 
a simply connected nilpotent Lie group and $\G$ is a cocompact lattice. 
Let $\mu$ be an ergodic invariant measure for $\a$. Suppose that the corresponding 
automorphisms of the Lie algebra $\g$ are diagonalizable, and either all their eigenvalues 
are real or $\mu$ is a weakly mixing measure for $\alpha$.

Then the conclusions (1) -- (4) of Theorem \ref{imNil} hold.
\end{corollary}
\smallskip

\begin{proof}
To apply Theorem \ref{imNil} we need to check for some regular element $b \in \Rk$ in the 
suspension of $\a$ that for each coarse Lyapunov subfoliation $W^P \subset W^- _b$ 
there exists a singular element $a$ in the corresponding Lyapunov hyperplane 
$\partial P$ which acts ergodically on $W^P$, i.e. $\xi _{a} \leq \xi (W^P)$.
In the case of TNS actions this is actually true for {\em any} regular element $b$
and any generic singular element $a\in \partial P$.

Consider a generic singular element $a\in \partial P$.
Take a regular element $c$ in the complement of $P$ so close to $a$ that 
it is not separated form $a$ by any Lyapunov hyperplane. Then all Lyapunov  
exponents that are nonzero for $a$ have the same signs on $a$ and $c$. 
Since the action is TNS, all Lyapunov exponents that correspond to $W^P$ become
positive on $c$. Thus $W_c^+=W_a^+\oplus W^P$ and $W_c^-=W_a^-$. Birkhoff averages 
with respect to $\a(a)$ of any continuous 
function are  constant on the leaves of $\xi_a$. Since such averages 
generate the algebra of $\a(a)$ invariant functions we conclude that 
$\xi_a\leq\xi(W_a^-)$. On the other hand both $\xi(W_c^-)$ and
$\xi(W_c^+)$ coincide with the Pinsker algebra $\pi(\a(c))$ (\cite{LY1}, 
Theorem B). Thus we conclude

$$\xi_a\leq\xi(W_a^-)=\xi(W_c^-)=\pi(\a(c))=\xi(W_c^+)=\xi(W_a^+\oplus
W^P)\leq\xi(E^P).$$ 
\end{proof}

\begin{corollary} \label{corKnil}
Suppose, in addition to the assumptions of Theorem \ref{imNil} or Corollary \ref{corTNS},
that $\a (c)$ is a K-automorphism with respect to the measure $\mu$ for some element 
$c \in \Zk$. Then $\mu$ is an algebraic measure, i.e. it is Haar on a single left 
coset $Hx$.  
\end{corollary}
\smallskip
We note that under K-automorphism assumption $\mu$ is automatically weakly mixing,
so that  complex eigenvalues are allowed with no extra assumption.

\begin{proof}
Since $\a (c)$ is a $K$-automorphism with respect to the measure $\mu$, the
factor with zero entropy has to be trivial. Hence $\mu$ is Haar on a single left 
coset $Hx$.
\end{proof}
\medskip

\subsection{Isomorphism Rigidity and Rigidity of Joinings}

As in the case of the homogeneous $\Rk$ actions, we apply the previous theorems
to obtain the following results on the isomorphism rigidity and rigidity of joinings.

\begin{theorem} \label{IrigNil}
Let $\alpha$ be a $\Zk$ action by automorphisms of $G/\Gamma$, where $G$ is a simply 
connected nilpotent Lie group and $\G \subset G$ is a cocompact lattice. Let $\lambda$ 
be the Haar measure on $G/\Gamma$. 
Let $(\alpha', \lambda ')$ be a $\Zk$-action of the same type on $G'/\G'$.

Suppose there is a subgroup $S \subset \Zk$ isomorphic to $\Z^2$ such that $\alpha (a)$ 
is ergodic with respect to $\lambda$ for all $a \in S$.     
Then any measure--preserving isomorphism $h: G/\Gamma \to G' / \Gamma '$  
between $(\a,\lambda)$ and $(\a', \lambda ')$ coincides \textup{(mod\;0)} with an 
affine isomorphism, and hence $\a$ and $\a '$ are algebraically isomorphic.
\end{theorem}

\begin{proof}
The proof of this theorem uses Corollary \ref{corKnil} and closely follows the proofs 
of Theorem \ref{Irig} and Corollary \ref{IrigSSimple}.  

Consider the suspension actions $\tilde \a$ and $\tilde \a '$ of $\Rk$
on the corresponding solvmanifolds $\M$ and $\M'$.
 By assumption, the restriction of $\a$ to $S \subset \Zk$
acts by ergodic automorphisms of $G/\G$. Then the proof of Corollary 4 in \cite{KS2} 
shows that the ergodicity assumption $\xi _{a} \leq \xi (W^P)$ of Theorem \ref{imNil}
is satisfied for the restriction of $\tilde \a$ to the linear hall $\tilde S \subset \Rk$ 
of $S$. We recall that any ergodic automorphism of $G/\G$ is a K automorphism
(\cite{P}). 

Consider the diagonal (product) action $\tilde \a \times \tilde \a'$ of $\R ^k$ 
on $\M \times \M'$. The measure--preserving isomorphism $h: G/\Gamma \to G' / \Gamma '$  
between $(\a,\lambda)$ and $(\a', \lambda ')$ can be lifted to a  measure--preserving 
isomorphism $\tilde h$ between $\tilde \a$ and $\tilde \a '$. Let $\tilde \mu$ be the push 
forward of the invariant measure on $\M$ to the graph of the measure--preserving isomorphism 
$\tilde h$. Then $\tilde \mu$ is an invariant measure for $\tilde \a \times \tilde \a'$ and 
can be viewed as the suspension of the invariant measure $\mu$ for $\a \times \a'$ which is 
the push forward of $\lambda$ to the graph of $h$ in $G/\Gamma \times G' /\Gamma '$. 

Since $(\a \times \a ', \tilde \mu)$ has the same ergodic properties as $\tilde \a$, 
the assumptions of Theorem \ref{imNil} and Corollary \ref{corKnil} are satisfied for 
the measure $\mu$ invariant under the restriction of $\a \times \a'$ to the subgroup 
$S \subset \Zk$. Thus $\mu$ is Haar on a single left coset $Hx \subset G/\Gamma \times 
G' /\Gamma '$. Since $h$ is a 
measure--preserving isomorphism and $\mu$ projects to the Haar measures in the factors, 
we see that the coset $Hx$ has to project one to one onto the factors. Thus $Hx$ is a 
graph of an affine isomorphism $h'$ between $G/\Gamma$ and $G' /\Gamma '$ which 
coincides \textup{(mod\;0)} with $h$.
\end{proof}
\medskip

\begin{theorem} \label{JrigNil}
Let $\alpha_1$ be a $\; \Zk$, $k \ge 2$, action by automorphisms of  $G_1/\Gamma_1$,
where $G_1$ is a simply connected nilpotent Lie group and $\G_1 \subset G_1$ is a 
cocompact lattice. Suppose that the corresponding automorphisms of the Lie algebra 
$\g_1$ are diagonalizable over $\C$ and suppose that there exists such an 
element $a\in \Zk$ that $\a_1 (a)$ is an ergodic automorphism with 
respect to the Haar measure $\lambda_1$ on $G_1/\G_1$.
Let $(\alpha_2,\lambda_2)$ be a $\Zk$-action on $G_2/\G_2$ of the same type.

Let $(G_1/\Gamma_1 \times G_2 /\Gamma _2, \a_1 \times \a_2, \mu)$ 
be a nontrivial joining of the actions $\a_1$ and $\a_2$. Suppose that either the 
eigenvalues of the corresponding automorphisms of $G_1/\Gamma_1 \times G_2 /\Gamma _2$ 
are all real,or $\mu$ is weakly mixing.

If for any one-parameter subgroup of $\Rk$, the ergodic components of the suspension of the 
joining measure $\mu$ contain  the fibers of the suspension, then the actions $\alpha_1$ and $\alpha_2$ have algebraic factors 
which  have algebraically isomorphic finite covers. Moreover, the joining measure $\mu$ is an extension 
of a zero-entropy measure for a common measurable factor with Haar measures in fibers. 
\end{theorem}

We prove this theorem in Section \ref{proof Jrig}.

\begin{corollary}[TNS Case] \label{JrigTNS}
Let $\alpha_1$ be a TNS $\; \Zk$ action by automorphisms of $G_1/\Gamma_1$,
where $G_1$ is a simply connected nilpotent Lie group and $\G_1 \subset G_1$ 
is a cocompact lattice. Suppose that the corresponding automorphisms of the 
Lie algebra $\g_1$ are diagonalizable over $\C$ and suppose that there exists 
 an element $a\in \Zk$ such that $\a_1 (a)$ is an ergodic automorphism with 
respect to the Haar measure $\lambda_1$ on $G_1/\G_1$.
Let $(\alpha_2,\lambda_2)$ be a $\Zk$-action on $G_2/\G_2$ of the same type.

Let $(G_1/\Gamma_1 \times G_2 /\Gamma _2, \a_1 \times \a_2, \mu)$ 
be a nontrivial joining of the actions $\a_1$ and $\a_2$. Suppose that either the 
eigenvalues of the corresponding automorphisms of $G_1/\Gamma_1 \times G_2 /\Gamma _2$ 
are all real, or $\mu$ is weakly mixing.
Then the actions $\alpha_1$ and $\alpha_2$ have algebraic factors which are 
algebraically isomorphic. Moreover, the joining measure $\mu$ is an extension 
of a zero-entropy measure for a common measurable factor with Haar measures 
in fibers.  
\end{corollary}

\begin{proof}
This corollary follows from Theorem \ref{JrigNil} similarly to the way 
Corollary \ref{corTNS} follows from Theorem \ref{imNil} (cf. \cite{KaK}).
\end{proof}

In the case of actions by automorphisms of tori the previous two results
can be slightly strengthened using the stronger theorem on invariant measures 
for the toral case, see \cite{KS4} Theorem 5.1' and \cite{KaK2} Theorem 3.1.  
The proof of Theorem \ref{JrigNil} yields in the toral
case Theorem \ref{ToralJoinGeneral} below. This theorem generalizes Theorem 3.1 
in \cite{KaK} and gives partial answer to the open problems posed in \cite{KaK}.
Instead of TNS assumption in Theorem \ref{ToralJoinGeneral} we may assume $k \ge 2$ 
and the ergodicity of one-parameter subgroups for the joining measure.

\begin{theorem} \label{ToralJoinGeneral}
Let $(\alpha_1,\lambda_1)$ and $(\alpha_2,\lambda_2)$
be TNS actions of $\Z^k$ by automorphisms
of $\T^{m_1}$ and $\T^{m_2}$, where $\lambda_1$ and $\lambda_2$ are
Lebesgue measures on $\T^{m_1}$ and $\T^{m_2}$ correspondingly. 
Suppose that there exists an element $a \in \Zk$ such that 
$\a_1(a)$ is ergodic with respect to $\lambda_1$.
If there exists a nontrivial joining measure $\mu$ on $\Tm =
\T^{m_1} \times \T^{m_2}$ then there exists a subgroup $\Gamma \subset \Zk$ 
of finite index such that the actions $\alpha_1$ and $\alpha_2$ restricted to 
$\Gamma$ have algebraic factors which are conjugate over $\Q$. 

Moreover, if the joining measure $\mu$ is ergodic it decomposes as 
$\frac1N (\mu_1 + ... + \mu_N)$, where each $\mu_i$, $i=1,...,N$,
is an invariant measure for the restriction of $\alpha_1 \times \alpha_2$ 
to $\Gamma$ and is an extension of a zero-entropy measure for the 
corresponding algebraic factor with Haar measures in the fibers. 
\end{theorem}
\bigskip

%%%%%%%%%%%%%%%%%%%%%%%%%%%%%%%% MAIN PROOF %%%%%%%%%%%%%%%%%%%%%%%%%%%%%%%%%%%

\section{Proof of Theorem \ref{im}} \label{proof im}

\subsection{Scheme of the Proof of Theorem \ref{im}} \label{im outline}
We will say that an integrable invariant distribution $E$ and its foliation $F$
satisfy property $\A$ if there exists a Lie subgroup $L \subset G$ such that 
for $\mu$ - a.e. $x$ the conditional measure $\m _x$ is a Haar measure on a single 
left coset $Lx$.

The main part of the proof is to show that $W_b^-$ satisfies property $\A$. This
is done in Section \ref{proof A}. After this we prove in Section \ref{constructionH} 
the existence of the subgroup $H$ which satisfies first three conclusions of the 
theorem. Then in Sections \ref{measurable factor} and  \ref{algebraic factor}
we describe the properties of the action in the factor and thus complete the proof
of the theorem. 

\subsection{The construction and properties of $H$} \label{constructionH}

Let $\L _0$ be the connected component of the identity of $\L (\mu)$, and let 
$M \subset \L _0$ be the maximal Lie subgroup of $\L _0$ generated by unipotent 
subgroups of $\L _0$. Since $W_b^-$ satisfies property $\A$ measure $\mu$ is
invariant under left translations by the corresponding Lie group $L$ and thus
$L$ is contained in $M$. Let $\mu _x$ denote the ergodic components of $\mu$ 
with respect to the left action of $M$.  
By Ratner's theorem  \cite[Corollary C]{R1} applied to $M$, each $\mu_x$ is 
algebraic, i.e. $\mu_x$ is Haar measure on some  closed orbit $H_x \cdot x$ 
of some subgroup $H_x$ of $G$ which contains  $M$. 

Since $A \subset \L _0$ and the 
adjoint action of $A$ maps unipotent elements to unipotent elements, we see
that $M$ is normalized by $A$. Hence the left action of $A$ preserves
the above ergodic decomposition, i.e.
                        $$a \mu_x = \mu _{a \, x}$$
In particular, $a$ maps the support of $\mu _x$ to the support of 
$\mu _{a \, x}$. Hence $a \, H_x  \, x = H_{a\,x}\, a\,x$ and thus
 $$H_{a\,x}\,= \,a\, H_x \, a^{-1}  \eqno (1)$$
We will now show that $H_x = H$ is $\mu$-a.e. constant. 
First we note that since $A$ acts ergodically (1) implies that $\dim(H_x)=d$ is 
$\mu$-a.e. constant. So we can consider the measurable function $\phi$ from $G/\Gamma$ 
to the Grassmanian $G_d$ of $d$-plains in $\g$ given by $x \mapsto \h_x$,
where $\h_x$ is the Lie algebra of $H_x$. By (1) we have 
$\phi (a\, x ) = a \, \phi (x) \, a^{-1}= Ad(a)\phi (x)$. Then $\nu=\phi_*(\mu)$ is
an ergodic measure on $G_d$ with respect to the adjoint action of $A$. To show that 
$H_x = H$ is $\mu$-a.e. constant we will prove that $\nu$ is a $\delta$-measure. 

Since for all $a \in A$, $Ad(a)$ is a diagonalizable over $\C$ automorphism of the 
Lie algebra $\g$, $\g$ splits into the direct sum of $\a$-invariant subspaces, 
which are the eigenspaces corresponding to real 
eigenvalues and invariant two-dimensional subspaces corresponding to the pairs of 
complex eigenvalues. Consider a regular element $c$ for which different Lyapunov
exponents take different values. The support of $\nu$ is contained in the 
set of non-wandering points of $Ad(c)$ on $G_d$. It is easy to see that this set 
consists of $d$-planes that are spanned by their intersections with the invariant 
subspaces of the above splitting. If eigenvalues of $Ad(a)$ are real for all $a \in A$
then clearly all such $d$-planes are fixed by the adjoint action of $A$. Hence the
ergodic measure $\nu$ which is supported on this set must be concentrated on a single
point. If complex eigenvalues are present, the set  of non-wandering points of $Ad(c)$
may include invariant tori on which the adjoint action of $A$ is represented by 
translations.
However, since in this case $\mu$ is assumed to be weakly mixing for $\alpha$ we see
that $\nu$ can not be supported on such tori.

Thus we conclude  that $H_x = H$ is $\mu$-a.e. constant. Now (1) implies
that $H_x = H_{a\, x} =  a \,H_x \, a^{-1}$ for $\mu$-a.e. $x$ and $a \in A$, 
i.e. $A$ normalizes $H$. This proves the first conclusion of Theorem~\ref{im}

Since the ergodic components of $\mu$ with respect to $M$ are Haar measures
on left cosets of $H$ we see that $\mu$ is invariant under left translations 
by $H$ and thus $H \subset \Lambda(\mu)$. This shows the second conclusion of 
Theorem~\ref{im}. 

Since $L$ is contained in $M \subset H$, its Lie algebra is 
contained in $E_b^- \cap \h$, where $\h$ is the Lie algebra of $H$.
But it clearly can not be strictly contained in $E_b^- \cap \h$, for 
otherwise the conditional measures on leaves of $W_b^-$ could not be supported
on a single left coset of $L$. This proves the third conclusion of Theorem~\ref{im}.

\subsection{Measurable factor with zero entropy.} \label{measurable factor} 

We have established above that for $\mu$- a.e. $x$ the ergodic component $\mu_x$ of 
$\mu$ with respect to $M$ is a Haar measure on the compact left coset $Hx$.
We note that this $\mu_x$ is also the ergodic component of $\mu$ with 
respect to $H$. Since $H$ is normalized by $A$ the left action of $A$ preserves
this ergodic decomposition. Hence  $\a$ induces a measurable factor 
action of $A$ on the space of the ergodic components of $\mu$ for the left action 
of $H$. Thus to establish the fourth  conclusion of Theorem~\ref{im} it remains to 
prove that the entropy of the factor measure is zero with respect to any element 
in the  Weyl chamber of $b$.

To prove that the entropy is zero we show that the measurable factor can be 
embedded into a topological factor of the restriction of $\a$ to the support of $\mu$
and then use Proposition 4.1 in \cite{KS3}.
\smallskip

We first show that the volume of cosets $Hx$ is constant almost everywhere.
The fact that the left action of $A$ preserves the above ergodic decomposition
means that 
                        $$a_* \mu_x = \mu _{a \, x}   \eqno (2)$$
where $\mu_x$ is the Haar measure on $Hx$ normalized to be probability. We fix some
(bi-invariant) Haar measure on $H$ and denote the corresponding measure on $Hx$
by $\nu_x$. Since both $\nu_x$ and $\mu_x$ are Haar measures $\nu_x=\rho (x)\mu_x$,
where $\rho (x)$ is an almost everywhere finite measurable function. Since
$a_* \nu_x = |det(Ad(a))| \nu _{a \, x}$ we obtain using (2) that 
$\rho (ax)= |det(Ad(a))| \rho (x)$. Using recurrence under $a$ we conclude that 
$|det(Ad(a))|=1$ $\mu$-a.e. This implies that $\rho (x)$ is $A$-invariant, hence
$Vol(Hx)=\nu_x(Hx)=\rho (x)$ is constant almost everywhere. We denote this constant
by $V$.

Now we show that for any $x \in supp (\mu) \subset \M=G/\G$ the coset $Hx$ 
has finite volume. For any $x \in \M$ the coset $Hx$, as an orbit of the left
action of $H$, is isomorphic to
$$H/{stab(x)} \cong H/{H \cap x\G x^{-1}}$$
Here we slightly abuse notations by writing $x$ in $x\G x^{-1}$ instead of 
an element in $G$ which projects to $x$. This is justified since $x\G x^{-1}$
does not depend on the choice of such an element. 

For any $x$ in the support of $\mu$ there exists a sequence of typical points
$x_n$ converging to $x$ for which $Vol(H/{H \cap x_n\G x_n^{-1}})=Vol(Hx_n)=V$. 
It is easy to see that the sequence of lattices $H \cap x_n\G x_n^{-1}$ converges 
to a discrete subgroup $\G'$ of $H$ which lies in $H \cap x\G x^{-1}$.
By Chabauty's Theorem \cite{Ragh}, 
$$Vol(Hx)=Vol(H/{H \cap x\G x^{-1}}) \le Vol(H/\G') \le$$
$$\liminf Vol(H/{H \cap x_n\G x_n^{-1}})=V.$$

Thus for any  $x$ in the support of $\mu$ the coset $Hx$ has volume bounded above 
by $V$. Now we show that it is compact. Since $G/\G$ is compact, there exists
a small open neighborhood $U$ of $e \in G$ such that 
$U \cdot U^{-1} \cap y\G y^{-1} =\{e\}$ for any $y \in G$. Hence $U \cdot U^{-1} \cap stab(y) =
U \cdot U^{-1} \cap H \cap y\G y^{-1} =\{e\}$. This implies that for any $y\in \M$ 
the projection from $H$ to $H/{H \cap y\G y^{-1}}$ is injective on $H \cap U$. 
Hence for any coset $Hx$ and point $y \in Hx$ there is a fixed size neighborhood 
of $y$ in $Hx$ which carries a fixed volume $Vol(H \cap U)$.

Since for all $x$ in the support of $\mu$ the coset $Hx$ has volume bounded above 
by $V$ we conclude that $Hx$ has to be compact. Moreover, the diameter of these
cosets (with respect to the intrinsic metric) is bounded uniformly in $x\in supp(\mu)$.

Now we can consider the projection $\pi$ given by $\pi(x)=Hx$, from the support of 
$\mu$ to the space $X$ of all compact subsets of $\M=G/\G$ equipped with the 
Hausdorff metric. Since the (intrinsic) diameters of these cosets are uniformly 
bounded in $x$ it is easy to see that $\pi$ is Lipschitz continuous.

Let us denote by $F$ and $G$ the projections to $X$ of the foliations $W_b^-$ and
$W_b^+ \oplus W_b^0$ respectively. Then $F$ is strictly contracted by any element $c$
in the Weyl chamber of $b$.  Since for $\mu$- a.e. $x\in \M$ 
the conditional measure of $\mu$ on the leaf $W_b^-(x)$ is Haar on a single left
coset $Lx$, where $L \subset H$, we conclude that the conditional measures of the
factor measure $\pi_*(\mu)$ on the leaves of $F$ are atomic almost everywhere.
Now Proposition 4.1 in \cite{KS3} implies that the entropy of the action in the
factor is zero for any element $c$ in the Weyl chamber of $b$. This completes
the proof of the fourth conclusion of Theorem~\ref{im}.

\subsection{Algebraic factor with zero entropy.}  \label{algebraic factor}

The last conclusion of Theorem~\ref{im} describes the case when the measurable
factor obtained above carries a natural algebraic structure. This is always
the case when the coset $N_G(H)y$ of the normalizer of $H$ in $G$ is compact
for some ``typical'' point $y$ whose $\a$-orbit is dense in $supp(\mu)$. 
Indeed, since $A$ normalizes $H$ the orbit $\a(y)$ is contained in  $N_G(H)y$. 
If the closure of the orbit $\a(y)$ is $supp(\mu)$ and if $N_G(H)y$ is compact
then clearly $supp(\mu) \subset N_G(H)y$.
In this case, the restriction of $\a$ to $N_G(H)y$ has an algebraic factor 
action isomorphic to $(y^{-1}Hy) \setminus y^{-1}N_G(H)y/ \G$. 
      
To complete the proof of the last conclusion of Theorem~\ref{im} it remains
to note that the entropy of the factor measure is zero with respect to {\it any} 
element in $A$. The fourth conclusion of Theorem~\ref{im} tells us that the entropy 
is zero for any element in the Weyl chamber of $b$. Clearly, the same is true for
any element in the Weyl chamber of $b^{-1}$. Since any element in $A$ can be
represented as a convex combination of elements in these Weyl chambers, the statement
follows from subadditivity of entropy for smooth actions established in \cite{H}. 

Since for any element $c \in A$ the entropy in the factor is zero, the conditional
measures on its stable foliation in this factor are atomic (\cite{KS3} Proposition 4.1) 
This implies that its stable foliation for the original action satisfies property $\A$.
Thus we have proved that the conclusion 3. of the theorem now holds for the stable 
foliation $W_c^-$ of {\bf any} element $c \in A$.

\bigskip
To complete the proof of Theorem~\ref{im} it remains to show that
$W_b^-$ satisfies property $\A$, which is done in the next section.

\bigskip

\section{$W_b^-$ satisfies property $\A$} \label{proof A}

The proof of the fact that $W_b^-$ satisfies property $\A$ has two main steps. 
First we show that each coarse 
Lyapunov subfoliation of $W_b^-$ satisfies property $\A$.
From this we deduce in  the second step that $W_b^-$ also satisfies property $\A$.

\subsection{\it Step 1} \label{step1}
In this section we show that each coarse 
Lyapunov subfoliation of $W_b^-$ satisfies property $\A$.
  
Let us consider one of the coarse Lyapunov subfoliation $W^P \subset W_b^-$,
where $P$ is the corresponding Lyapunov half-space that contains $b$.
By the assumption of Theorem \ref{im} there exists a singular element $a$ in the
Lyapunov hyperplane $\partial P$ which acts ergodically on $W^P$, i.e.
$\xi _{a_i} \leq \xi (W^P).$
Recall that we denote by $\xi_a$ the partition into ergodic components
of element $a$ and by $\xi(F)$ the measurable hull of $F$.

The desired result is established by applying Lemmas \ref{cm1}, 
\ref{cm2}, \ref{cm3}, and \ref{cm4} to $F=W^P$. Note that by the assumption of 
Theorem \ref{im}, $Ad(c)$ is a diagonalizable automorphism of the Lie algebra $\g$ 
for all $c \in A$. This implies that the derivative of the element $a$ acts
isometrically on $W^P$.  

The general setup for these lemmas is as follows. 
Let $a \in A$ be a singular element acting isometrically on an invariant
foliation $F$. 
Denote by $B _1^F (x)$ the closed unit ball about $x$ in the leaf $F(x)$
with respect to the induced metric. Denote by $\m _x$ the
system of conditional measures on $F$ normalized by the requirement
$\mu _x (B^F_1 (x) ) =1$ for all $x$ in the support of $\mu$. Denote by
$I_x$ the subgroup of all isometries of $F(x)$ which preserve $\m _x$ up to 
a scalar multiple.

The first three lemmas are the analogs of the Lemmas
5.4, 5.5, and 5.6 from \cite{KS3}, see also Lemmas 3.2, 3.3, and 3.4
in \cite{KaK2}.

\begin{lemma} \label{cm1} (\cite{KS3}, Lemma 5.4) Suppose that $a$ acts isometrically 
on some $\a(a)$-invariant foliation $F$ which satisfies $\xi_a \leq \xi (F)$. 
Then for $\mu$-a.e. $x$, $I_x$ is closed and the support of $\m _x$ is the orbit of 
$x$ under the group $I_x$. Furthermore, $ \phi_* \m _x = \m _{\phi x}$ for any $\phi \in I_x$.
\end{lemma}

The proof of Lemma 5.4 in \cite{KS3} (Lemma 3.2 in \cite{KaK2}) works
verbatim for the general case.

\begin{lemma}   \label{cm2} Suppose that an invariant subfoliation $F$ of $W_b^-$ 
for an element $b$ satisfies the conclusions of Lemma \ref{cm1}. Then for $\mu$-a.e. 
$x$, the support $S_x$ of $\m _x$ is a homogeneous submanifold of $F(x)$, 
more precisely it is a left coset $S_x=L_x x$ of some subgroup $L_x$.  
\end{lemma}

\begin{proof} 
Let $E$ denote the tangent distribution of the foliation $F$. Let us denote by 
$\lambda_1, ..., \lambda_m$ the negative Lyapunov exponents of $b$ and by
$E^{\lambda_i}$, $i=1,...,m$, the intersections of $E$ with the corresponding 
Lyapunov subspaces.

The conclusion of Lemma \ref{cm1} states that for $\mu$-a.e. $x$ the support
$S_x$ of the measure $\m _x$ is the orbit of a closed
group of isometries. Therefore $S_x$ is a submanifold of the leaf $F(x)$,  
possibly disconnected. Thus the tangent space $T_x S_x$ depends on $x$
smoothly along the leaf $F(x)$ and measurably on $\M$. Our main goal it to show 
that the tangent distribution of $S_x$ is invariant under left translations along 
$S_x$ in the following sense. If $x$ and $hx$ belong to the same connected component 
then the left translation by $h$ takes the tangent space at $x$ to the tangent space 
at $hx$.

We note that $b$ maps $S_x$ to $S_{bx}$ and contracts $F$. Using recurrence under 
the iterates of $b$ and the fact that the tangent vectors in different $E^{\lambda_i}$
have different exponential rates of contraction one can see that 
$T_x S_x= \bigoplus_{i=1}^m (T_x S_x \cap E^{\lambda_i})$. Thus to prove that the tangent 
distribution $T S_x$ is invariant under left translations it suffices to show that
$T_x S_x \cap E^{\lambda_i}$ is invariant under left translations for every $i=1,...,m$.

Let us fix $i$ and put $V_x=T_x S_x \cap E^{\lambda_i}$. Note that for $y \in S_x$ 
we have $S_x=S_y$. Let $\beta_x (y)$ be the maximal angle between the subspaces $V_y$ 
and $R_g (V_x)$, where $g=x^{-1}y$. Then
$$K(x)= \limsup_{\epsilon \to \, 0} \frac1\epsilon 
\sup \{\beta_x (y) \, | \,  d(x,y) < \epsilon \}$$
gives a maximal rate of deviation of the distribution $V_x$ from the left 
invariant distribution 
and thus plays the role of maximal curvature. Since $V_x=T_x S_x \cap E^{\lambda_i}$, 
it is measurable in $x$ and smooth along $S_x$. Thus we see that
$\beta_x (y)$ is measurable in $x$ and Lipschitz along $S_x$. Thus $K(x)$
exists and is measurable. Moreover, it is easy to see that $K(y)$ gives
an upper bound for the derivative of $\beta_x$ at $y$. Hence to prove that 
$V_x$ is invariant under left translations for almost every $x$ it suffices 
to show that $K(x)=0$ almost everywhere.

To show this we use the element $b$ for which $\lambda (b) <0$. Thus the
iterates of $b$ exponentially contract the leaves of $F$. We note that $b$ 
maps $V_x$ to $V_{b\,x}$ and left invariant distributions to left invariant 
distributions. Since $V_x$ is contained in a single Lyapunov subspace,  
we see that the exponential contractions in all directions $V_x$ are 
the same. Since the derivative of $b$ along $E^{\lambda_i}$ is semi-simple, 
$b$ preserves the angles. Since $b$ also contracts the distances it follows 
that $K(b ^n x)$ goes to infinity unless $K (x)=0$. 
Since Poincar\'{e} recurrence under $b$ implies that this is impossible, we
conclude that $K(x)=0$ almost everywhere. Hence $V_x$ is left invariant.
\bigskip

Thus we conclude that the tangent distribution of every connected component
of $S_x$ is left invariant in the sense that if $x$ and $hx$ belong to the 
same connected component then the left translation by $h$ takes the tangent
space at $x$ to the tangent space at $hx$. 
Hence $S_x= \bigcup x_i H_{x_i}$.

Let us now show that the support is connected. Suppose to the
contrary that the support  is a union $\cup A_i$ of at least two
affine subspaces $A_i$.  Let $d_x$ denote the minimum of the distances
from $x$ to any $A_i$ which does not contain $x$.  Since the support
is a closed subset, $d_x >0$ for all $x$. Note that $d_{b^n \, x}
\rightarrow 0$ as $n \rightarrow \infty$. This again contradicts
Poincar\'{e} recurrence under $b$.  

\end{proof}

\noindent {\bf Remark}. 
If we do not assume that the derivative of $b$ along $E^{\lambda_i}$ is semi-simple, 
it may contain Jordan blocks. In this case the angles may not be preserved, but the 
distortion of the angles is at most polynomial and does not affect the above argument 
since the distances are contracted exponentially. 
\bigskip

\begin{lemma}       \label{cm3} Suppose that a $b$-invariant subfoliation 
$F \subset W_b^-$ satisfies the conclusions of Lemma \ref{cm1} and suppose that  
the support $S_x$ of $\m _x$ is, for $\mu$-a.e.$x$, a left coset $S_x=L_x x$ 
of some subgroup $L_x$. Then $\m _x$ is Haar measure on $S_x=L_x x$.
\end{lemma}

A modified version (see Lemma 3.2 in \cite{KaK2}) of the proof of 
Lemma 5.4 in \cite{KS3} can be applied for the general case.
It shows that the isometries of Lemma \ref{cm1} actually preserve the measure.
We note that the existence of such a transitive group of measure preserving
isometries implies that the measure is Haar.

\begin{lemma}       \label{cm4}  Suppose that for an $\a$-invariant foliation $F$ 
the conditional measure $\m _x$ is, for $\mu$-a.e. $x$, Haar measure on a left 
coset $S_x=L_x x$ of some subgroup $L_x$. Then $L_x=L$ is constant almost everywhere.
Hence $F$ satisfies property $\A$ and $aLa^{-1}=L$. 
\end{lemma}      

\begin{proof}
The proof is the same as the proof in Section \ref{constructionH} of the fact that 
$H_x$ is constant almost everywhere.
\end{proof}

\bigskip

\subsection{\it Step 2} \label{step2}
In this step we complete the proof of Theorem \ref{im} by showing that the foliation
$W_b^-$ satisfies property $\A$.
We will prove this by using the results of Step 1 and establishing property $\A$
inductively for larger and lager invariant subfoliations of $W_b^-$.

Let $P_i$, $1 \leq i \leq l$, be the Lyapunov half-spaces containing $b$ and
$E^{P_i}$ ($W^{P_i}$) be the corresponding coarse Lyapunov distributions in
$E_b^-$ (subfoliations of $W_b^-$). The results of the first step imply that 
all foliations $W^{P_i}$ satisfy property $\A$.

To arrange the inductive process we restrict the action to a generic 2-plane which 
contains $b$ and intersects all the Lyapunov hyperplanes in generic lines.
For each Lyapunov hyperplane $\partial P_i$ we choose a generic singular element
$c_i$ in the intersection of $\partial P_i$ with the generic 2-plane.
We can now reorder elements $c_i$ and Lyapunov half-spaces $P_i$ in such a way that 
for any $1 \leq j \leq l-1$ the distribution $E^{P_1} \oplus ... \oplus E^{P_j}$ is 
integrable and is contracted by $c_{j+1}$.

Now the fact that the foliation $W_b^-$ satisfies property $\A$ follows
from inductive application of the following lemma, which is the main part of Step 2. 

\begin{lemma}    \label{bigstep} 
Let $W$ be an invariant subfoliation of $W_c^I$ and $F$ be an invariant subfoliation 
of $W_c^-$ for some generic singular element $c$. Suppose that $W$ and $F$ satisfy 
property $\A$ and 
are jointly integrable to the foliation $W \oplus F \subset W_b^-$ for some element $b$. 
Then $W \oplus F$ satisfies property $\A$
\end{lemma}

\begin{proof}
We prove this lemma inductively by adding one dimensional subfoliations of $W$
to the foliation $F$ until we exhaust the whole $F \oplus W$. We first set up
the inductive process and then use Lemma~\ref{smallstep} to show that on each 
step of this process we obtain a foliation that satisfies property $\A$.

{\bf Notation.}  Throughout the proof we denote by $H$ the one-dimension\-al 
subfoliation of $W$ which is being added at the current step of the induction and 
by $D$ the sum of all one-dimensional subfoliations of $W$ which have been
already added (thus $D$ is trivial in the beginning). We denote by $E^F$, $E^W$, 
$E^D$, and $E^H$ the tangent distributions of foliations $F$, $W$, $D$, and $H$
correspondingly.
\smallskip

{\bf Description of the induction.} 
To arrange the induction we will choose a basis of $E^W$ and an ordering of
this basis according to which the corresponding one-dimensional subfoliations 
will be added. 
Since $W\subset W_b^-$ is an invariant subfoliation, $E^W$ is an invariant 
subalgebra of $E_b^-$. Since $E_b^-$ is nilpotent so is $E^W$. Thus the
central series of $E^W$ gives the invariant filtration $E^W=V^k \supset ...  
\supset V^1 \supset V^0=0$ such that $[V^i,V^i] \subset V^{i-1}$, for $i=1,...k$. 
In particular each $V^i$ is a subalgebra invariant under the adjoint action of
$b$ and $c$. Recall that $L(\mu)$ is the Lie algebra of $\L(\mu)$, the subgroup of $G$
whose left translations preserve $\mu$. Hence $V^i \cap L(\mu)$ is a subalgebra
of $V^i$, for $i=1,...k$, which is invariant since $L(\mu)$ is. To obtain the
desired ordered basis of $E^W$ we first choose a basis 
of $V^1 \cap L(\mu)$ and complement it to a basis of $V^1$. After that
we complement it to a basis of $V^1 \oplus V^2 \cap L(\mu)$ and then
to a basis of $V^1 \oplus V^2$. Continuing in the same way we obtain a basis
of $E^W$ with the natural partial order. The order of the basis directions within 
$V^1 \cap L(\mu)$ and  within each step of the complementing can be chosen
arbitrarily. We note that the subspace spanned by the first several basis 
directions is always a subalgebra and thus is integrable. Hence the foliations $D$ 
and $D \oplus H$ above are well defined.  

Since $Ad(b)$ and $Ad(c)$ are diagonalizable, commute, and leave all subspaces  
$V^i$ and $V^i \cap L(\mu)$ invariant, the basis above can be chosen in such 
a way that on each inductive step the foliations $D$ and $D \oplus H$ are  
invariant for some multiples $sb$ and $tc$ in $A$ of the elements $b$ and $c$. 
Indeed, if all eigenvalues of $Ad(b)$ and $Ad(c)$ on $E^W$ are real the basis 
above can be chosen as a common eigenbasis for $Ad(b)$ and $Ad(c)$. In this case
on each step the foliations $D$ and $D \oplus H$ are clearly invariant under
$b$ and $c$. If $Ad(b)$ and $Ad(c)$ have complex eigenvalues, a similar basis 
can be chosen using the real normal forms for $Ad(b)$ and $Ad(c)$. We then add 
the one-dimensional foliations in such a way that if one direction of a complex 
eigenspace gets added, the other gets added in the next step. Now choose 
multiples $sb$ and $tc$ of $b$ and $c$ in such a way that their eigenvalues 
corresponding to $E^H$ are real. Then on each inductive step the foliations 
$D$ and $D \oplus H$ are invariant under  $sb$ and $tc$.

Thus we have arranged the inductive process in such a way that on each step the 
foliations  $D$ and $D \oplus H$ are well defined and invariant for some multiples 
$sb$ and $tc$ of the elements $b$ and $c$. We also note that the distributions 
$E^F \oplus E^D$, $E^F \oplus E^H$ and $E^F \oplus E^D \oplus E^H$ are integrable. This clearly 
follows from the next observation. If a subdistribution $E$ of $E^W$ is a subalgebra then 
the direct sum $E \oplus E^F$ is also a subalgebra and thus is integrable. Indeed, 
since both $E$ and $E^F$ are subalgebras it suffices to show that the Lie bracket 
of an element in $E$ with an element in $E^F$ belongs to $E^F$. Since $F \oplus W$ 
is integrable it must belong to $E^F \oplus E^W$. We note that if the Lie bracket 
is nonzero its Lyapunov exponent is the sum of the Lyapunov exponents of the original 
elements. Since the kernels of these two Lyapunov exponents are different we see that 
the kernel of their sum is different from both of them. Since $E^W \subset E^I_c$ and
all Lyapunov exponents of $E^I_c$ have the same kernel we see that the Lyapunov subspace 
of the bracket cannot be in $E^W$. Hence it must belong to $E^F$.
\smallskip

{\bf Inductive step.} 
It follows from the construction of the basis that on each step either $E^H$ 
is contained in $L (\mu)$ or the intersection of $E^D \oplus E^H$ with $L (\mu)$ 
is contained in $E^D$. 

If the current $H$ is contained in $L (\mu)$ the inductive step is simple.  
Indeed, in this case for $\mu$-a.e. $x$ the conditional measure on the leaf 
$F \oplus D \oplus H(x)$ is invariant under the left translations
by the Lie subgroup $U \subset G$ corresponding to $(E^{F \oplus D} \cap L(\mu)) 
\oplus E^H$. Hence this conditional measure is an integral of Haar measures on
the left cosets of $U$. Moreover, we see that it must be supported on a single left 
coset of $U$ since we know that the conditional measure on the leaf 
$(F \oplus D) (x)$ is supported on a single left coset of $U \cap (F \oplus D)$.
Thus $F \oplus D \oplus H$ satisfies property $\A$.

The next lemma proves the inductive step in the case when the intersection of 
$E^D \oplus E^H$ with $L (\mu)$ is contained in $E^D$. This will complete
the proof of Lemma~\ref{bigstep}.
\end{proof}

\begin{lemma}  \label{smallstep} Let $F$, $D$, and $H$ be $b$ and $c$ invariant
subfoliations of $W_b^-$. Suppose that $H$ is one-dimensional, that the distributions 
$E^F \oplus E^D \oplus E^H$, $E^F \oplus E^D$, $E^F \oplus E^H$ and $E^D \oplus E^H$ are integrable, 
and that we have $(D \oplus H) \subset W_c^I$ and $F \subset W_c^-$.

If the foliations $F \oplus D$ and $D \oplus H$ satisfy property $\A$ and if
the intersection of $E^D \oplus E^H$ with $L (\mu)$ is contained in $E^D$, then
the foliation $F \oplus D \oplus H$ also satisfies property $\A$.  
\end{lemma}

\begin{proof} 
We consider a measurable partition, which is subordinate to the foliation 
$F \oplus D \oplus H$ and tiles the intersection of each leaf with a set of measure 
at least $0.99$ by small ``boxes'' of roughly the same size. The boxes can be chosen
in the following way. Take $y \in \M$ and small open neighborhoods $U$ and $V$ of $y$
in $F(y)$ and $D(y)$ correspondingly. Since $F \oplus D$ is integrable we can define
a rectangle $R_y$ as the set of intersections   
$$\{z \in (F \oplus D)(y)\, :\; z=F(v) \cap D(u), \text{ where } 
v \in V \text{ and } u \in U\}.$$
Similarly, since $(F \oplus D) \oplus H$ is integrable we can define a box $T_y$ 
as the product of the above type of $R_y$ with an interval in one-dimensional leaf
$H(y)$. For any $x$ in the above set of measure at least $0.99$ above we denote by $T(x)$
the box it belongs to. 

Since the foliation $D\oplus H$  satisfies property $\A$ there exists
such a compact set $K$ of measure at least $0.98$ that for any point $x$ 
the intersection of the box $T(x)$ with $K \cap (D\oplus H)(x)$ 
is contained in the intersection of $T(x)$ with a single left coset of 
$\L (\mu) \cap (D \oplus H)$. Since the foliation $F \oplus D$ also 
satisfies property $\A$ we may assume the similar relationship between $K$
and $F \oplus D$.

Denote be $X$ the set of points in $K$ whose distance from the union of the boundaries
of the boxes is at least $\gamma$. By the boundary of a box we understand its leafwise
boundary with respect to the leaf of $F \oplus D \oplus H$. The measurable partition 
can be easily chosen in such a way that the total measure of the union of these boundaries 
is small and $\mu (X) > 0.96$ for some $\gamma >0$.

Let $\mu _X$ be the restriction of the measure $\mu$ to the ``good'' set $X$, i.e.
$\mu _X (\cdot)=\mu (X \cap \cdot)$.
We note that the measure $\mu _X$ may not be invariant. Let us consider the system 
of the conditional measures of $\mu_X$ with respect to the measurable partition into 
the boxes described above (the set not covered by the boxes has $\mu_X$ measure 0). 
These measures will be referred to as the conditional measures of boxes. 

We try to slice each box $T(x)$ into three parts $T_x^L$, $T_x^M$, and  $T_x^R$ 
of equal conditional measure by two leaves of the foliation $F \oplus D$, here
$T_x^M$ is the middle part. This may not be possible if the intersection
of a box with a single leaf of $F \oplus D$ has positive conditional measure,
but then we can proceed as in the second case of the dichotomy below.
For each $x \in X$ we denote by $d(x)$ the 
minimal distance (along the one-dimensional $H$ direction) between $T_x^L$ and  
$T_x^R$. For the measurable function $d(x)$ we have the following dichotomy:

\begin{enumerate}

\item Either $d(x) \ge d>0$ on a set $Y$ which consists of whole 
boxes and has measure $\mu_X (Y) > 0.95$,

\item or for every box in a set of positive measure, at least 1/3 of its 
conditional measure is concentrated on the intersection with a single leaf 
of $F\oplus D$.

\end{enumerate}

In the latter case we argue as follows. Consider a point $x$ in this set of 
positive measure and the corresponding box $T(x)$.
Since the foliation $F \oplus D$ 
satisfies property $\A$ the conditional measure on the leaf $(F \oplus D)(x)$
is Haar on a single left coset of $\L (\mu) \cap (F \oplus D)$. The fact that
the intersection of this left coset with the box $T(x)$ carries at
least 1/3 of the conditional measure of $T(x)$ means that the measure of 
this intersection dominates the measure of its small tubular neighborhood
in the leaf $(F \oplus D \oplus H)(x)$. By applying measure preserving 
transformation $b^{-1}$, which expands $F \oplus D \oplus H$, we see that, 
in a set of positive measure, the Haar measure on a single left coset of 
$\L (\mu) \cap (F \oplus D)$ dominates larger and larger part of the leaf of 
$F \oplus D \oplus H$. Using this and recurrence under $b^{-1}$ we conclude
that the conditional measure on a typical leaf of $F \oplus D \oplus H$ is
Haar on a single left coset of $\L (\mu) \cap (F \oplus D)$. This proves
that $F \oplus D \oplus H$ satisfies property $\A$. 

To complete the proof it remains to eliminate the first alternative of the 
dichotomy. Now we use the special properties of $D \oplus H$.  

Recall that the intersection of the compact set $K$ with any box 
consists of the pieces of left cosets of $\L (\mu) \cap D$,
with at most one piece for each leaf of $D\oplus H$. Hence the projection (along
$D$ direction) of this intersection to the ``$F\oplus H$ face of the box''
is the graph of a continuous (not necessarily everywhere defined) function 
from the  $F$-direction to the $H$-direction. Since $K$ is compact the family of these 
functions is equicontinuous.  We will show that this equicontinuity 
contradicts the recurrence under the action of a properly chosen element in $A$.

Since $b$ contracts $F \oplus D \oplus H$ we can find a sufficiently
large number $n$ such that the size of the image of any box
under the action of $b^n$ is $\gamma$-small.  The distance along
the $H$ direction between the images of $T_x^L$ and  
$T_x^R$, which was at least $d$, becomes at least $d'$. Let
us fix some $\epsilon<d'$ and consider $\delta>0$ given for this $\epsilon$
by the equicontinuity. Since $c$ acts isometrically on $D \oplus H$
and contracts $F$ we can choose $k$ so large that the image of any
box under $c^k b^n$ is $\delta$-small in $F$ direction. We conclude
that the element $c^k b^n$ satisfies the following conditions:

\begin{enumerate}
\item The image of any box is $\delta$-narrow in the $F$
direction

\item The distance along the $H$ direction between the images of
the right and the left parts of any box is at least $\epsilon$

\item The diameter of the image of any box is less than $\gamma$.  
\end{enumerate} 

Consider a point $y \in Y$ and the corresponding box $T$. Since the set $X$
does not intersect the $\gamma$-neighborhood of the boundaries of the boxes,
3. implies that $c^kb^n(T) \cap X$ lies in one box. Hence by 1. and 2., $c^kb^n(T^L)$ and
$c^kb^n(T^R)$ cannot both intersect $X$. This implies that there exists 
a set $Z \subset Y$ with $\mu_X(Z) \ge \frac13 \mu_X(Y)$ for which
$c^kb^n(Z) \cap X= \emptyset$. But this is impossible since $\mu(X)>0.96$
and $\mu(c^kb^n(Z)) = \mu(Z) \ge \mu_X(Z) \ge \frac13 \mu_X(Y) > \frac13 0.95$.
 
This shows that the second alternative of the dichotomy is impossible and thus 
proves that foliation $F \oplus D \oplus H$ satisfies property $\A$.

\end{proof}

This completes  the proof of Theorem \ref{im}.

%%%%%%%%%%%%%%%%%%%%%% OTHER PROOFS %%%%%%%%%%%%%%%%%%%%%%%%%%%%%%%%

\section{Proof of Corollary \ref{corSSimple}.} \label{proof corSSimple}

We will use the notations from Theorem \ref{im} and its proof. 

We first claim that the subgroup $H$ is reductive. W.l.o.g. we suppose
that $H$ itself rather than a conjugate intersects $\Gamma$ in a lattice. 
If $H$ is not reductive, then it can be decomposed as $H= S \, K \, Rad$ where $S$ is
semisimple without compact factors, $K$ is a compact semisimple group, and
$Rad$ is the the radical of $H$. Then $K \, Rad$ intersects $\Gamma$ in a 
lattice (cf. e.g. \cite[Corollary 2.12]{W}). That $H$ is reductive follows 
once we know that the unipotent radical $R_u$ of $H$ is trivial. If not, we 
may pick a nontrivial 
one parameter unipotent subgroup $U$ of $R_u$. Any nontrivial $ u \in U$
is horocyclic in the sense of \cite[Definition 2]{R3}, as follows
from the Jacobson-Morozov lemma. Since $\Gamma$ is a cocompact lattice in
$G$, we may apply \cite[Corollary 2]{R} to the subgroup $K\,
Rad$. Hence $K \, Rad$ will contain a copy of $sl_2$ which is impossible.
Thus $H$ is reductive. 

Recall that the Haar measure on $H \, x$ is ergodic w.r.t. a subgroup
generated by unipotents contained in $H$. In particular, it is ergodic
w.r.t. $S$. Note that $S$ acts trivially on the maximal Euclidean quotient
of $H \, x$ (the quotient by $S \, K$). Hence $S$ cannot act ergodically
unless $H$ is semisimple. 

Since $H$ is semisimple the connected component of its normalizer is a
product of $H$ with the centralizer of $H$ (connected components always).
Thus 
$H$ intersects $\Gamma$ in a lattice which is Zariski dense in $H$  by the
Borel density theorem. Hence the centralizer of $H$ is the same as the
centralizer of
$H \cap \Gamma$. The latter is an intersection of centralizers of all the
elements
of $H \cap \Gamma$. Since $\Gamma$ is arithmetic by Margulis' arithmeticity
theorem, $G$ is defined over a number field $F$, and
$\Gamma$ commensurable with the integer points of $G$ w.r.t. $F$. Hence the
centralizer of $\gamma
\in \Gamma \cap H$ is defined over $F$, and so is $Z_G H$. Hence $\Gamma
\cap N_G H$ is commensurable with integer points of $N_G H$ which form a
lattice on $N_G H$ by  Borel's and Harish Chandra's theorem \cite{BHC}.

\bigskip

\section{Proof of Theorem \ref{imNil}}  \label{proof imNil}

In this section we will reduce Theorem  \ref{imNil} to  Theorem \ref{im}.
As $\alpha$ is given by diagonalizable automorphisms, the homomorphism from $\Zk$ to the automorphisms of $G$
extends to $\Rk$. Hence we can form the semidirect product groups $G' = \Rk \ltimes G$, $\Gamma ' =  \Zk \ltimes G$
and its compact quotient $G' /\Gamma '$. Then the   action $\alpha '$ of $\Rk$ on $G' /\Gamma '$ is the  
action of $\Rk$ induced from $\Zk$.  We will freely use all the notation from
the main theorem, denoting any relevant object for the induced action by a ' . 

Then $\alpha '$ almost satisfies the assumptions of Theorem \ref{im}.
The main assumption on the existence of suitable singular elements is already postulated in the formulation of
Theorem \ref{imNil}. If  the eigenvalues of the automorphisms $\alpha (\Zk)$ have real eigenvalues so do those of
 $\alpha (\Rk)$. Thus we only need to discuss what to do when   $\alpha $ is weakly mixing.
 Consider the natural  fibration of $G' / \Gamma '$ over $T= \Rk / \Zk$. 
Then any eigenfunction $f$ for $\Rk$ is an eigenfunction for $\Zk$ on a.e. fiber $G/\Gamma$. 
Hence $f$ is constant along a.e. fiber. In the proof of the main theorem, we consider the 
function $x \mapsto H_x '$. If $\alpha$ is weakly mixing, the $H_x '$ are constant 
along a.e. fiber $G/\Gamma$, as desired.  Since the $H_x '$ are generated
by unipotent
elements, it is clear that $ H_x ' \subset G$. 

Note that conclusions 1.-3. of Theorem \ref{im} immediately imply the corresponding claims of Theorem \ref{imNil}.
Furthermore, statements 4. and 5. of Theorem \ref{im} imply conclusion 4. of Theorem \ref{imNil} if we can show that
the normalizer $N_{G'} (H')$ intersect $\Gamma '$ in a lattice. This follows from the following result. We thank 
 D. Witte for providing a proof for us. As we could not find precise references for standard
 facts, we  give a detailed account.

\begin{proposition} 
 Suppose
 \begin{itemize}
 \item $G$ is a $1$-connected, nilpotent Lie group,
 \item $A$ is a $1$-connected, abelian Lie group which acts on $G$ by automorphisms, 
 \item $G' = A \ltimes G$,
 \item $H$ is a closed, connected subgroup of $G$,
 \item $\Gamma $ is a lattice in $G$,
 \item $\Gamma \cap H$ is a lattice in $H$,
 \item $\Zk$ is a lattice in $A$.
 \end{itemize}
Let $\Gamma ' = \Zk \ltimes \Gamma$.  Then $\Gamma ' \cap N_{G'}(H)$ is a lattice in $N_G(H)$,
and also $\Gamma  \cap N_G(H)$ is a lattice in $G$. 
\end{proposition}

\begin{proof}

Since $\Gamma '$ is a polycyclic group, $\Gamma ' $ admits a faithful representation $\rho $ into 
$GL(n,\Z)$ by \cite[Theorem 4.12]{Ragh} such that the image of the nilradical of $\Gamma '$ consists
of unipotent matrices. Thus $\rho (\Gamma) $ is unipotent. As $\Gamma$ is nilpotent, $\rho$ extends
to a representation of $G$ into a unipotent subgroup of $GL(n,\R)$  by \cite[Theorem 2.11]{Ragh}. 
Also $\rho \mid \Zk$ extends to a representation of $A$ into $GL(n,\R)$ such that $\rho (A)$
normalizes $\rho (G)$. 

Hence we may assume $G' \subset GL(\ell,\R)$, for some $\ell$, in such a way
that
 \begin{itemize}
 \item $G$ is unipotent, and
 \item $G ' \cap GL(\ell,\Z)$ is (commensurable with) $\Gamma$.
 \end{itemize}
 Let 
 \begin{itemize}
 \item $L$ be the Zariski closure of $G'$ in $GL(\ell,\R)$, and
 \item $\Lambda = L \cap GL(\ell,\Z)$.
 \end{itemize}
  Dani's generalization of the  Borel Density Theorem \cite{Dani2}  shows
   that $\Lambda$ is Zariski dense in a
cocompact subgroup of $L$, which implies that $L$ has no characters
defined over $\Q$. By the compactness criterion for lattices \cite[Theorem 10.19]{Ragh}, 
$\Lambda$ is a cocompact lattice in $L$.

Because $\Gamma \cap H$ is Zariski dense in $H$, we know that $H$ is
defined over $\Q$. Hence $N_L(H)$ is also defined over $\Q$.
Therefore, $\Lambda \cap N_L(H)$ is a cocompact lattice in $N_L(H)$.
By \cite[Lemma 1.19]{Ragh}, we conclude that $\Lambda \cap G ' \cap N_L(H)$ is a
lattice in 
 $G' \cap N_L(H) = N_ 
 {G'} (H)$. This proves the first claim since  $\Lambda \cap G '$ is commensurable
with $\Gamma '$. Finally, that $\Gamma  \cap N_G(H)$ is a lattice in $G$ follows immediately
since $N_G(H) / \Gamma  \cap N_G(H)$ is the fiber of the fibration $N_{G'} (H) / \Gamma  \cap N_{G'}
(H)$ over the torus $\Rk / \Zk$. 

\end{proof}

\section{Proof of Theorems \ref{JrigSSimple} and \ref{JrigNil}.}  \label{proof Jrig}
 
In this section we prove Theorems \ref{JrigSSimple} and \ref{JrigNil}.
The proofs are very similar, so we first prove Theorem \ref{JrigSSimple}
and then in Section \ref{proof JrigNil} describe the changes needed for
the nilmanifold case.

\subsection{Proof of Theorem \ref{JrigSSimple}.}  \label{proof JrigSSimple}
We introduce the following notations: $\a=\a_1 \times \a_2$, $G=G_1 \times G_2$, 
$\G=\Gamma_1 \times \Gamma_2 \subset G$, and $A=\rho(\Rk)=\{(\rho_1(a),\rho_2(a)) 
\in G : a \in \Rk \}$, where $\rho =\rho_1 \times \rho_2$ is the diagonal embedding 
of $\Rk$ into $G=G_1 \times G_2$.

Since $(\a, \mu)$ satisfies the assumptions of Corollary \ref{corSSimple}
there exists a subgroup $H \subset \Lambda (\mu) \subset G$ such that $\mu$
is supported on a compact coset $N_G(H)y \subset G/\G$ of the normalizer of 
$H$ in $G$. By passing if necessary to the algebraically isomorphic action of 
$y^{-1}Ay$ on $(y^{-1}Hy) \setminus y^{-1}N_G(H)y/ \G$, we may assume that 
$\bar e =\{e\G\}\in G/\G$ is a density point for $\mu$, and hence $\mu$
is supported on $N_G(H) \bar e \subset G/\G$. 

Let $H_i =H \cap G_i$, $i=1,2$. 
We would like to factorize by $H_1 \times H_2 \subset H$, so we first prove that
$H_1 \times H_2$ is normal in $G$. By symmetry, it suffices to show that $H_1$ is 
normal in $G_1$. For this we note that since $\mu$ is supported on $N_G(H)$ and 
projects to the Haar measure on the $G_1/\G_1$, the projection of $N_G(H)$ to $G_1$
must be onto. Hence for any $g_1 \in G_1$ there exists  $g_2 \in G_2$ such that
$(g_1,g_2)\in N_G(H)$, i.e. $(g_1^{-1},g_2^{-1})H(g_1,g_2)=H$. This, together with
the equality $(g_1^{-1},g_2^{-1})G_1(g_1,g_2)=G_1$, shows that 
$(g_1^{-1},g_2^{-1})(G_1 \cap H)(g_1,g_2)=(G_1 \cap H)$. But the latter implies that
$g_1^{-1} H_1 g_1= H_1$. Thus we have proved that $H_1 \times H_2$ is normal in $G$.

Hence we can consider the factor action $\a'$ on $G'/\G'$, where 
$G'=(H_1 \times H_2)\setminus G=G_1' \times G_2'$ with $G_i'=H_i\setminus G_i$, $i=1,2$. 
We note that since $\G =\G_1 \times \G_2$, the lattice $\G'$ in the factor is also
the product of lattices $\G_1'$ and $\G_2'$, which are the projections of $\G_1$ and 
$\G_2$ to $G_1'$ and $G_2'$ respectively.
Recall that for the measure $\mu$ we have a natural measurable partition $\xi$ into 
compact cosets $Hx$, $x \in N_G(H)$, with Haar conditional measures. It follows 
that the factor measure $\mu'$ on $G'/\G'$ is invariant under the left action of 
$H'=(H_1 \times H_2)\setminus H$, and that the projection of $\xi$ to $G'/\G'$ is a 
measurable partition $\xi'$ into compact cosets $H'x$. We note that since 
$H' \cap G_i' =\{e\}$, 
the subgroup $H'$ projects injectively to $G_1'$ and $G_2'$. Let $H_1'=p_1(H')$,
where $p_1$ is the projection to the first factor. Since $p_1(H'x)=H_1' \cdot p_1(x)$,
it is easy to see that $\xi'$ projects to a measurable partition $\xi_1'$ of 
$\M_1'=G_1'/\G_1'$ which consists of compact cosets $H_1'y$. We will 
show below that $\xi_1'$ is coarser than the Pinsker partition for $\a_1'(a)$, 
where $\a_1'$ is the action on $\M_1'$ and $a \in \Rk$ is from the 
statement of Theorem \ref{JrigSSimple}. Since $\a_1'(a)$ is a factor of $\a_1(a)$ and the
latter is, by assumption, a K automorphism with respect to the Haar measure on 
$G_1/ \G_1$, we conclude
that $\xi_1'$ has to be the trivial partition. Therefore, $\mu_1'=(p_1)_* \mu'$ is
concentrated on a single compact coset $H_1'x$. Since $\mu_1'$ is Haar, as the 
projection to $\M_1'$ of the Haar measure $\lambda_1$ on $G_1 / \G_1$, we conclude that
$H_1'=G_1'$. Thus the subgroup $H'$ projects injectively and onto $G_1'$. By the
same argument, $H'$ projects injectively and onto $G_2'$. Hence $H'$ is the graph
of the isomorphism $f : G_1' \to G_2'$ determined by the equality 
$H'=\{(g,f(g))\in G': g\in G_1'\}$.

We know that $\mu'$ is invariant under $H'$ and is supported on the normalizer
$N_{G'}H'$ of $H'$ in $G'$. We now describe $N_{G'}H'$ to show that $\mu'$ is supported 
on a single coset $H'x$. If $(g_1,g_2) \in N_{G'}H'$, then for any $g \in G_1'$ we
have $(g,f(g))\in H'$ and thus $(g_1^{-1}gg_1,g_2^{-1}f(g)g_2)$ is again in $N_{G'}H'$. 
Hence $f(g_1^{-1}gg_1)=g_2^{-1}f(g)g_2)$ which implies 
$f(g)=(f(g_1)g_2^{-1})^{-1} f(g) (f(g_1)g_2^{-1})$. Since this holds for all $g \in G_1$,
$z=(f(g_1)g_2^{-1})$ lies in the centralizer of $f(G_1')$, which coincides with the center
$Z(G_2')$ of $G_2'$ as $f$ is an isomorphism. Since $z\cdot (g_1,g_2)=(g_1,f(g_1))\in H'$
we conclude that $N_{G'}H'=Z(G_2)\cdot H'$. Since $G_2'$ is semisimple, $Z(G_2')$ is
at most countable and thus $N_{G'}H'$ consists of at most countably many cosets of $H'$.
Since $\mu'$ is an ergodic measure for the $\Rk$ action $\a'$ supported on  $N_{G'}H'$,
it is clear that $\mu'$ is supported on a single coset $H'x \subset G'/\G'$. 
Since $\bar e \in G'/\G'$ is a density point of $\mu'$, we see that $\mu'$ is Haar on 
$H'\bar e$. This implies that $\mu$ is Haar on $H \bar e$, i.e. $\mu$ is algebraic. 

Since $\mu'$ is Haar on $H'\bar e$, the latter is also invariant under $\a'$. Since
we know that $H'$ is isomorphic to $G_1'$ and $G_2'$ it follows that $H' \bar e$ projects 
onto $\M_1'$ and $\M_2'=G_2'/(\G' \cap G_2'$. 
In contrast to the case of isomorphism rigidity $H'x$ may not project injectively to 
$\M_1'$ and $\M_2'$. Thus it may not produce an algebraic isomorphism between the actions 
$\a_1'$ and $\a_2'$. However, $H'x$ is compact and thus $H' \cap \G'$ is a uniform 
lattice in $H'$. Hence the projection of $H' \cap \G'$ to $G_1'$ is also a uniform 
lattice which is contained in $\G_1'$ since $\G' =\G_1' \times \G_2'$, and thus is
of finite index in $\G_1'$. Similarly, the projection of $H' \cap \G'$ to $G_2'$ is
of finite index in $\G_2'$. This implies that the invariant coset $H'\bar e$ is a 
graph of an algebraic isomorphism between the actions $\a_1'$ and $\a_2'$ on the 
finite coverings of $G_1'/\G_1'$ and $G_2'/\G_2'$ respectively.
\medskip

To complete the proof it remains to show that the measurable partition $\xi_1'$ of 
$\M_1'$ is coarser than the Pinsker partition for $\a_1'(a)$. We note that the Pinsker 
partition for $\a_1'(a)$ coincides with the measurable hull of the unstable foliation
for $\a_1'(a)$ (\cite{LY1} Theorem B). Thus it suffices to prove that the projection
$p_1(\h')$ of the Lie algebra $\h'$ of the group $H'$ to the Lie algebra $\g_1'$ 
of the group $G_1'$ contains the Lie algebra of the unstable distribution of
$\a_1'(a)$. Consider a Lyapunov exponent $\chi$ for the action $\a'$ such that
$\chi(a)>0$ and denote by $V^\chi$ the corresponding Lyapunov subspace in the
Lie algebra $\g'$ of the group $G'$. Then $\chi$ is also a Lyapunov exponent for
for at least one of the factors, and $V^\chi$ is the direct sum of the corresponding
Lyapunov subspaces $V^\chi_1 \subset \g_1'$ and $V^\chi_2 \subset \g_2'$, one of which
may be trivial. We recall that both projections $p_1$ and $p_2$ are injective on $H'$.
Hence 
$$\dim (\h' \cap V^\chi) \le \dim V^\chi_i, \quad i=1,2. \qquad \eqno (1)$$  

Since $\mu_1'$ is the Haar measure on $\M_1'$, the entropy of $\a_1'(a)$
is $$h(\a_1'(a))= \sum \chi(a) \dim (V^\chi_1),  \qquad \eqno (2)$$
where the sum is taken over all Lyapunov exponents $\chi$ positive on $a$.
Since for $\mu$-a.e. $x$ the conditional
measure of $\mu$ on the leaf $W_{\a(a)}^+(x)$ is Haar on a single left coset $Hx$, it follows
that for $\mu'$-a.e. $x$ the conditional measure of $\mu'$ on the leaf $W_{\a'(a)}^+(x)$ is 
Haar on a single left coset $H'x$. This implies (\cite{LY2}) that the entropy of $\a'(a)$ is 
$$h(\a'(a))= \sum \chi (a) \dim (H \cap V^\chi),  \qquad  \eqno (3)$$  
where the sum again is taken over all Lyapunov exponents positive on $a$.
Since $\a_1'(a)$ is a factor of $\a'(a)$, $h(\a_1'(a)) \le h(\a'(a))$. 
Thus (1), (2), and (3) imply that $\dim (\h' \cap V^\chi)= \dim V^\chi_1$ for
all Lyapunov exponents $\chi$ positive on $a$. Since $p_1$ is injective on 
$\h' \cap V^\chi$, this shows that it is also surjective. Hence $p_1(\h')$ contains
the Lyapunov subspaces of all Lyapunov exponents positive on $a$ and thus contains
the unstable distribution $E^+_{\a_1'(a)}$. This completes the proof of 
Theorem \ref{JrigSSimple}.

\subsection{Proof of Theorem \ref{JrigNil}.}  \label{proof JrigNil}

In this section we show how to adapt the proof of Theorem \ref{JrigSSimple} above
to the case of $\Zk$ actions by automorphisms of a nilmanifold, and thus
establish Theorem \ref{JrigNil}.

Let $\a=\a_1 \times \a_2$, be the diagonal action of $\Zk$ on $G_1/\G_1 \times 
G_2/\G_2$, i.e. $\a(a)(x_1,x_2)=(\a_1(a)x_1,\a_2(a)x_2)$. Then the joining 
$(\a, \mu)$ satisfies the assumptions of Theorem \ref{imNil}. Hence there 
exists an $\a$ invariant subgroup $H \subset \Lambda (\mu) \subset G_1 \times G_2$ 
such that $\mu$ is supported on a compact coset $N(H)y$ of the normalizer of $H$ 
in $G_1 \times G_2$. Again by passing if necessary to the algebraically isomorphic 
action of $\Zk$ on $(y^{-1}Hy) \setminus y^{-1}N_G(H)y/ \G$, we may assume that 
$\mu$ is supported on $N(H) \bar e$, where $\bar e$ is the projection of the unit
to  $G_1/\G_1 \times G_2/\G_2$. We again can factorize by $H_1 \times H_2 \subset H$,
where $H_i =H \cap G_i$, $i=1,2$. In the same way we obtain that the measure $\mu'$
in the factor is Haar on the compact cosets of $H'=(H_1 \times H_2) \setminus H$,
and that $H'$  projects isomorphically to the factors $G_1'$ and $G_2'$ of $G_1$ 
and $G_2$ respectively. We recall that the ergodicity of the element $\a_1(a)$
from the statement of Theorem \ref{JrigNil} implies that it is a K-automorphism
(\cite{P}).

In contrast to the semisimple case above, in this case $\mu'$ 
may not be supported on the single coset, so that $\mu'$ and $\mu$ may not be algebraic.
According to Theorem \ref{imNil} the measures $\mu'$ and $\mu$ are extensions of
a zero entropy measure for a measurable factor of $\a$ with Haar conditional measures
in the fibers. This measurable factor can be viewed as an
algebraic factor of the restriction of $\a$ to $N(H) \bar e$. Even though $\mu'$ 
may not be supported on a single coset of $H'$, we already know that $H$ and hence 
$H'$ are $\a$ invariant, and the coset $H' \bar e$ is compact. As in the proof above, 
these facts easily imply that $H' \bar e$ is a graph of an algebraic isomorphism between 
the factor actions on the finite coverings of $G_1'/\G_1'$ and $G_2'/\G_2'$, where 
$\G_i'$ is the projection of the lattice $\G_i$ to $G_i'$, $i=1,2$. The measurable
factor above can be viewed as a common measurable factor of these factor actions.

%%%%%%%%%%%%%%%%%%%%%%%%%%%%%%%%%%%%%% END %%%%%%%%%%%%%%%%%%%%%%%%%%%%%%%%%%%%%%%%%%

\end{document}